\documentclass{amsart}

\usepackage{amssymb,amsfonts, latexsym,amsmath,hhline,array,longtable}
\usepackage{pdfsync,color, comment,colortbl, mathrsfs,stmaryrd,cite,graphicx}

\usepackage{ifpdf}
\ifpdf \usepackage[colorlinks=true, citecolor=blue, linkcolor=blue, urlcolor=blue]{hyperref} \fi

\newcommand{\cal}{\mathcal}

\newtheorem{formula}{}[section]
\newtheorem{definition}[formula]{Definition}
\newtheorem{corollary}[formula]{Corollary}
\newtheorem{remark}[formula]{Remark}
\newtheorem{lemma}[formula]{Lemma}
\newtheorem{theorem}[formula]{Theorem}

\def\thrm{\begin{theorem}}
\def\thrml#1{\begin{theorem}\label{#1}}
\def\ethrm{\end{theorem}}
\def\rmrk{\begin{remark}}
\def\rmrkl#1{\begin{remark}\label{#1}}
\def\ermrk{\end{remark}}
\def\dfntn{\begin{definition}}
\def\dfntnl#1{\begin{definition}\label{#1}}
\def\edfntn{\end{definition}}
\def\nmrt{\begin{enumerate}}
\def\enmrt{\end{enumerate}}
\def\tm#1{\item[{\rm (#1)}]}
\def\qtnl#1{\begin{equation}\label{#1}}
\def\eqtn{\end{equation}}
\def\lmm{\begin{lemma}}
\def\lmml#1{\begin{lemma}\label{#1}}
\def\elmm{\end{lemma}}
\def\crllr{\begin{corollary}}
\def\crllrl#1{\begin{corollary}\label{#1}}
\def\ecrllr{\end{corollary}}
\def\css{\begin{cases}}
\def\ecss{\end{cases}}
\def\prf{\begin{proof}}
\def\eprf{\end{proof}}

\def\cD{{\cal D}}

\def\cL{{\cal L}}

\def\cX{{\cal X}}
\def\cY{{\cal Y}}

\def\mF{{\mathbb F}}
\def\mN{{\mathbb N}}

\def\fF{{\mathfrak F}}
\def\fK{{\mathfrak K}}

\def\fX{{\mathfrak X}}

\DeclareMathOperator{\aut}{Aut}

\DeclareMathOperator{\GQ}{GQ}

\DeclareMathOperator{\id}{id}
\DeclareMathOperator{\im}{Im}

\DeclareMathOperator{\mon}{Mon}

\DeclareMathOperator{\pr}{pr}

\DeclareMathOperator{\sym}{Sym}

\DeclareMathOperator{\WL}{WL}

\def\mltst#1{\{\hspace{-3pt}\{#1\}\hspace{-3pt}\}}

\def\qaq{\quad\text{and}\quad}
\def\qoq{\quad\text{or}\quad}
\def\ov{\overline}

\begin{document}

\title[A large family of strongly regular graphs]{A large family of strongly regular graphs with small Weisfeiler-Leman dimension}

\author{Jinzhuan Cai}
\address{Hainan University, Haikou, China}
\email{caijzh12@163.com}
\author{Jin Guo}
\address{Hainan University, Haikou, China}
\email{guojinecho@163.com}
\author{Alexander L. Gavrilyuk}
\address{Shimane University, Matsue, Japan}
\email{gavrilyuk@riko.shimane-u.ac.jp}
\author{Ilia Ponomarenko}
\address{Hainan University, Haikou, China; Steklov Institute of Mathematics at St. Petersburg, Russia}
\email{inp@pdmi.ras.ru}

\thanks{}
\date{}


\begin{abstract}
In 2002, D. Fon-Der-Flaass constructed a prolific family of strongly regular graphs. In this paper, we prove that for infinitely many 
natural numbers $n$, this family contains $n^{\mathsf{\Omega}(n^{2/3})}$ strongly 
regular $n$-vertex graphs~$X$ with the same parameters, which satisfy 
the following condition: an isomorphism between $X$ and any other 
graph can be verified by the $4$-dimensional Weisfeiler-Leman algorithm.
\end{abstract}

\maketitle

\section{Introduction}
An undirected graph $X$ is said to be \emph{strongly regular} if the number $k$ (respectively, $\lambda$, $\mu$) of common neighbors of any two of its vertices depends only on whether they are equal (respectively, adjacent, nonadjacent). Together with the order of $X$, these numbers form the set of \emph{parameters} of~$X$ and are intersection numbers of the coherent configuration associated with~$X$ (the exact definitions are given in Section~\ref{240623a}). 
To date, the monograph \cite{BrouwerM2022} can serve as 
the most complete (though not exhaustive) reference for 
the theory of strongly regular graphs.

Informally speaking, ``strongly regular graphs lie on the cusp 
between highly structured and unstructured'' \cite{Cameron}. 
Given the parameters $(n,k,\lambda,\mu)$, 
there may be very few (or none at all) pairwise nonisomorphic 
strongly regular graphs having these parameters or, conversely, 
an abundance of those, say, the number of them is bounded from 
below by an exponential function in~$n$.

In the former case, a strongly regular graph $X$ can typically be
nicely characterized in terms of its local structure, which gives 
rise to polynomial-time algorithms to test isomorphism between~$X$ 
and any other graph (examples here are the Johnson, Hamming, and 
Grassman graphs of diameter two). Moreover, the role 
of such an algorithm is often played by a low-dimensional 
Weisfeiler-Leman algorithm \cite{EvdP2000,Ponomarenko2020a}; 
in other words, the Weisfeiler-Leman dimension 
(see Section~\ref{240623b}) of the graph $X$ is usually bounded 
from above by a small constant.

In the latter case, there are numerous examples, such as 
the graphs of Latin squares, Steiner systems, and the graphs 
obtained by various prolific constructions  \cite{FonDF2002,Ihringer2017,Muzychuk2009F,VVK,BIK}.
However, the authors were not aware of any abundant family 
of strongly regular graphs with small Weisfeiler-Leman dimensions. 
This was a primary motivation for the  paper and our main result is to present such a family. 

\thrml{050623a}
For any large enough power of $2$, say $q$, 
there is a family $\fF=\fF_q$ of pairwise nonisomorphic 
strongly regular graphs with the same parameters,  
such that 
\qtnl{220623g}
\dim_{\scriptscriptstyle\WL}(\fF)\le 4\qaq|\fF|>q^{\mathsf{\Omega}(q^2)},
\eqtn
where $\dim_{\scriptscriptstyle\WL}(\fF)$
is the maximum Weisfeiler-Leman dimension over the graphs 
in~$\fF$.
\ethrm

The parameters of all $n$-vertex  strongly regular graphs in  the family $\fF$ from Theorem~\ref{050623a} are equal to  
\qtnl{220623a}
(n,k,\lambda,\mu)=(q^2(q+2),q(q+1),q,q ).
\eqtn 

In fact, $\fF$ is a part of a larger prolific  family 
of  strongly regular graphs constructed by D.~Fon-Der-Flaass in~\cite{FonDF2002}. Recall that his construction has the following three degrees of freedom:
\nmrt
\item[$\bullet$] one can choose arbitrary affine designs $A_1,\ldots,A_{q+2}$ with the same parameters, each~$A_i$ has  $q+1$ parallel classes, 
\item[$\bullet$] one can choose an arbitrary coloring of the arcs  of a directed complete graph $K_{q+2}$  in $q+1$ colors so that the arcs coming from any vertex have different colors,
\item[$\bullet$] for any two disinct indices $i$ and $j$, one can 
choose a bijection between the lines of the parallel classes 
specified by a pair~$A_i,A_j$.
\enmrt

For the family of graphs constructed in Theorem~\ref{050623a}, 
each degree of freedom is reduced significantly. First, all the designs~$A_i$ are equal to the same Desarguesian affine plane $A$ of even order $q$ ($q>16$ is sufficient for the assumption of Theorem~\ref{050623a}).
Second, the coloring of the arcs of the complete graph $K_{q+2}$ is defined with the help of an arbitrarily chosen hyperoval\footnote{A subset of~$A$ consisting  of $q+2$ points  such that no three pairwise distinct points belong to the same line.} in~$A$; the details are given in Section~\ref{230623j}. Finally, the bijections between the lines 
of the parallel classes are not arbitrary. 
In fact, all the graphs in $\fF$ are obtained from one distinguished strongly regular graph $X^*$ by switching of a special type, see Section~\ref{290623a}. Though $X^*$ does not belong to the family~$\fF$, it has the same parameters as the graphs in $\fF$ and belongs to the Fon-Der-Flaass family; moreover, $X^*$ is the collinearity graph of 
a generalized quadrangle $\GQ(q+1,q-1)$, 
see~\cite[Example~10.18] {KissS2020}. 

Our result 
leaves open several natural questions. 
Denote by $\fF’$ the family of all Fon-Der-Flaass graphs constructed from a fixed affine plane~$A$ and a fixed arc coloring of the graph $K_{q+2}$. Then the second part of  Eq. \eqref{220623g} and the calculations given in \cite[Proposition~3.5]{Muzychuk2009F} show that the quotient $\frac{|\fF|}{|\fF’|}$ tends to zero as the number~$q$ increases. 
We wonder if it is possible to show, perhaps using a more general type of switching, that the Weisfeiler-Leman dimension of almost all graphs of the family  $\fF'$ is bounded from above by some absolute constant.

In connection with this question, it is worth noting that the Fon-Der-Flaass graphs constructed with $A_1,\ldots,A_{q+2}$ being the same affine design with even number of points and the high-dimensional graphs 
of Cai-F\"urer-Immerman type \cite{Fuhlbr2021}, 
whose Weisfeiler-Leman dimension is unbounded, 
are special cases of the same general  construction of graphs 
(not necessarily strongly regular), described in~\cite[Section~3]{Muzychuk2009F}. In this construction, the number of copies of 
the design does not exceed the number $a$ of its parallel classes, increased by~$1$. Strongly regular graphs are obtained when $a=q+1$, and the graphs from \cite{Fuhlbr2021}  are obtained when the number~$a$ coincides with the dimension of the affine space. In at least these two examples, the properties of a graph to be strongly regular or to have a large Weisfeiler-Leman dimension lie, so to speak, on the opposite ends of the spectrum.

It would also be interesting to see whether there exists a larger family of strongly regular Fon-Der-Flaass graphs under some milder 
restrictions than the ones we have chosen, e.g., if the affine designs are not necessarily affine planes, or if they are planes but not necessarily Desarguesian. Note also the Fon-Der-Flaass construction  
was generalized in several ways, see \cite{Muzychuk2009F,VVK}.
\smallskip

The most challenging step in the proof of Theorem~\ref{050623a} is to verify the first inequality in Eq.~\eqref{220623g}. Here we use the theory of coherent configurations (see Section~\ref{240623a}). Within the framework of this theory, every graph $X$ is associated with a  coherent configuration $\cX=\WL(X)$ constructed by the $2$-dimensional Weisfeiler-Leman algorithm. Moreover, the Weisfeiler-Leman dimension of the graph $X$ can be estimated from above by the base number $b(\cX)$ of the coherent configuration $\cX$, increased by two, where the \emph{base number} is the smallest number of points with respect to which the extension of~$\cX$ is a discrete coherent configuration (in fact, we prove a slightly stronger analog of this result, see Lemma~\ref{090623a}). The 
switching method, which we used to construct graphs of the family $\fF$, allows us to estimate $b(\cX)$ from above by the base number of the coherent configuration associated with the affine plane~$A$, which, as is well known, does not exceed~$2$, see \cite[Theorem~3.3.8]{CP2019}.~
\smallskip

We conclude the introduction by expressing 
our deepest gratitude to M.~Muzychuk, who not only drew our attention to the study of the Fon-Der-Flaass graphs, but also explained numerous subtle details of the construction.

\section{Graphs and coherent configurations}\label{240623a}

\subsection{Notation.}
Throughout the paper, $\Omega$ stands for a finite set. For any $\Delta\subseteq \Omega$, we denote by~$1_\Delta$ the diagonal of the Cartesian square $\Delta\times\Delta$, and abbreviate $1_\alpha:=1_{\{\alpha\}}$ for $\alpha\in\Delta$. The set of all classes
of an equivalence relation $e$  on a subset of $\Omega$ is  denoted by $\Omega/e$.

For a binary relation $r\subseteq\Omega\times\Omega$, we set 
$r^*=\{(\beta,\alpha)\colon\ (\alpha,\beta)\in r\}$. The set of all {\it neighbors} of a point $\alpha\in\Omega$ in the relation $r$ is denoted by $\alpha r=\{\beta\in\Omega\colon\ (\alpha,\beta)\in r\}$, 
and the number $|\alpha r|$ is the \emph{valency} of $\alpha$ in $r$. 
For relations $r,s\subseteq\Omega\times\Omega$, we put 
$$
r\cdot s=\{(\alpha,\beta)\colon\ (\alpha,\gamma)\in r,\ (\gamma,\beta)\in s\text{ for some }\gamma\in\Omega\},
$$ 
which is called the \emph{dot product} of $r$ and $s$.
For $\Delta,\Gamma\subseteq \Omega$, we set $r_{\Delta,\Gamma}=r\cap (\Delta\times \Gamma)$ (and abbreviate $r_{\Delta}:=r_{\Delta,\Delta}$) and 
$\Delta r=\cup_{\delta\in \Delta}\delta r$. 
For a set $S$ of relations on $\Omega$, 
we denote by $S^\cup$ the set of all unions of the elements of $S$ and 
put $S^*=\{s^*\colon s\in S\}$ and $S^f=\{s^f\colon s\in S\}$ for any bijection $f$ from $\Omega$ 
to another set.

Finally, $\mathsf{\Omega}()$ stands for 
the Big-Omega (a tight lower bound) notation.

\subsection{Graphs}
By a {\it graph} we mean a (finite) simple undirected graph, i.e., a pair $X=(\Omega,E)$, where $E\subseteq \Omega\times \Omega$ is reflexive and symmetric. The elements of  the sets~$\Omega$ and~$E$ are called, respectively, the {\it vertices} and {\it edges} of the graph~$X$. Two vertices $\alpha$ and $\beta$ are said to be {\it adjacent} (in $X$) whenever $(\alpha,\beta)\in E$ (equivalently, $(\beta,\alpha)\in E$). The subgraph of $X$ induced by a set $\Delta\subseteq \Omega$ is denoted by $X_{\Delta}=(\Delta,E_{\Delta})$.

The set $\alpha E\cap\beta E$ of all common neighbors of the vertices $\alpha$ and $\beta$ in the graph $X$ is denoted by~$N_X(\alpha,\beta)$. The graph $X$ is said to be {\it strongly regular} with parameters $(n,k,\lambda,\mu)$ if $|\Omega|=n$ and  the number
$$
n_X(\alpha,\beta)=|N_X(\alpha,\beta)|
$$
is equal to $k$, $\lambda$,  or $\mu$ depending on whether the vertices $\alpha$ and $\beta$ are equal, adjacent, or non-adjacent.  

Following  \cite[Section 4]{HestH1971}, a graph $X$ satisfies the {\it $4$-condition} if   the number of $4$-vertex subgraphs of~$X$ of a given type with respect to a given pair $(\alpha,\beta)$ of distinct vertices depends only on whether $\alpha$ and $\beta$ are adjacent or not; here, two subgraphs of $X$ are of the same type with respect to  a  pair $(\alpha,\beta)$ of distinct vertices if both contain $\alpha$ and $\beta$ and there exists an isomorphism from one onto the other 
mapping $\alpha$ to $\alpha$ and $\beta$ to $\beta$. 
In any strongly regular graph satisfying the $4$-condition, 
the number
$$
e_X(\alpha,\beta)=|E_{N_X(\alpha,\beta)}|
$$
depends only on whether distinct vertices $\alpha$ and $\beta$ are adjacent or not.

\subsection{Coherent configurations}
Let $S$ be a partition of the Cartesian square $\Omega^2$; in particular, the elements of $S$ are treated as binary relations on~$\Omega$. A pair $\cX=(\Omega,S)$ is called a {\it coherent configuration} on $\Omega$ if the following conditions are satisfied:
\nmrt
\tm{C1}  $1_{\Omega}\in S^\cup$,
\tm{C2} $S^*=S$,
\tm{C3} given $r,s,t\in S$, the number $c_{rs}^t=|\alpha r\cap \beta s^{*}|$ does not depend on $(\alpha,\beta)\in t$. 
\enmrt

Any relation belonging to $S$ (respectively, $S^\cup$)  is called a {\it basis relation} (respectively, a {\it relation} of~$\cX$). The set of all relations is closed with respect to taking the transitive closure and the dot product. A set $\Delta \subseteq \Omega$ is called a {\it fiber} of $\cX$ if the relation $1_{\Delta}$ is basis. The set of all fibers is denoted by $F=F(\cX)$. Any element of $F^\cup$ is called a {\it homogeneity set} of $\cX$. From \cite[Theorem~2.6.7]{CP2019}, it follows that given $s\in S^\cup$ and $d\in \mN$, we have
\qtnl{180723d}
\{\alpha\colon|\alpha s|\leq d\}\in F^\cup.
\eqtn 

There is a partial order\, $\le$\, of the coherent configurations on the same set~$\Omega$. Namely, given two such coherent configurations $\cX$ and $\cX'$, we set $\cX\le\cX'$ if and only if each basis relation of~$\cX$ is the union of some basis relations of~$\cX'$. The minimal and maximal elements with respect to this ordering are the {\it trivial} and {\it discrete} coherent configurations: the basis 
relations of the former one are the reflexive relation $1_\Omega$ and its complement 
in $\Omega\times\Omega$ (if $|\Omega|\geq 1$), whereas the basis relations of 
the latter one are singletons. The discrete coherent configuration on $\Omega$ is denoted by $\cD_\Omega$.

Given an affine plane $A$ of order $q$, one can define a coherent 
configuration $(\Omega,S)$ called the \emph{affine scheme}\footnote{See  \cite[Subsection~2.5.2]{CP2019}.} 
of the plane. It has degree $|\Omega|=q^2$ and rank $|S|=q+2$. 
The points of $\Omega$ are just the points of $A$, while 
the irreflexive basis relations in $S$ correspond to the parallel 
classes $A$. Namely, every parallel class $P$ defines the set of 
pairs $(\alpha,\beta)$ (with $\alpha\ne \beta$) such that 
the line through $\alpha,\beta$ belongs to $P$.

Let $\cX=(\Omega,S)$ and $\cX'=(\Omega',S')$ be two coherent configurations. A bijection $f\colon\Omega\to\Omega'$ is called an {\it isomorphism} from $\cX$ to $\cX'$ if $S^f= S'$. The isomorphism $f$ induces a natural bijection $\varphi\colon S\to S'$, $s\mapsto s^f$. One can see that $\varphi$  preserves the numbers from the condition~(C3), namely, the numbers $c_{rs}^t$ and $c_{r^{\varphi},s^{\varphi}}^{t^{\varphi}}$ are equal  for all $r,s,t\in S$. Every bijection $\varphi\colon S\to S'$ having this property is called an {\it algebraic isomorphism}, written as $\varphi\colon \cX\to \cX'$. The algebraic isomorphism $\varphi\colon\cX\to\cX'$ induces a uniquely determined bijection $S^\cup\to {S'}^\cup$ denoted also by $\varphi$.

A coherent configuration is called {\it separable} if every algebraic isomorphism from it to another coherent configuration is induced by an isomorphism. The  trivial and discrete coherent configurations are separable, see \cite[Example~2.3.31]{CP2019}.

The {\it coherent closure} $\WL(r,s,\ldots)$ of the  binary relations $r,s,\ldots$ on $\Omega$, is defined to be the smallest\footnote{with respect to the natural	partial order on the partitions of the same set.} 
coherent configuration on $\Omega$, containing each of them as a relation. When $r=E$ is the edge set of a graph $X$, we write $\WL(X,s,\ldots):=\WL(r,s,\ldots)$.  The {\it coherent configuration of a graph} $X$ is just the coherent closure of its edge set: $\WL(X)=\WL(E)$. Note that the coherent closure  is a closure operator on the set of all partitions of $\Omega^2$ satisfying conditions (C1) and (C2).

The \emph{extension} $\cX_{\alpha_1,\ldots,\alpha_b}$ of a coherent configuration $\cX$ with respect to points $\alpha_1,\ldots,\alpha_b$ is defined to be $\WL(X,1_{\alpha_1},\ldots,1_{\alpha_b})$. The \emph{base number} $b(\cX)$ of a coherent configuration $\cX$ is 
the minimal integer $b\ge 0$ such that the extension of $\cX$ 
with respect to some $b$ points is discrete.


\section{The Weisfeiler-Leman dimension of a graph}\label{240623b}
Throughout this section, $m\ge 2$ is an integer and  $M=\{1,\ldots,m\}$.  The elements of the Cartesian $m$-power $\Omega^m$ are represented by the $m$-tuples $x=(x_1,\ldots,x_m)$ with $x_i\in\Omega$ for all $i\in M$.  

\subsection{Multidimensional Weisfeiler-Leman algorithm}
The Weisfeiler-Leman dimension of a graph is defined with the help of 	the $m$-dimensional Weisfeiler-Leman algorithm. For a  given graph $X=(\Omega,E)$, this algorithm  constructs a certain coloring $c(m,X)$ of  $\Omega^m$; here, a \emph{coloring} is understood 
as a function $c$ from ~$\Omega^m$ to some linearly ordered 
set $\im(c)$ whose elements are called \emph{colors}.

At the first stage, an initial coloring $c_0=c_0(m,X)$ is constructed  from an auxiliary coloring~$c'$ defined as follows: given $x,y\in\Omega^m$,  we set $c'(x)=c'(y)$ if and only if for all $i,j\in M$,
\qtnl{250623k}
x_i=x_j\ \Leftrightarrow\ y_i=y_j\qaq  (x_i,x_j)\in E\ \Leftrightarrow\ (y_i,y_j)\in E.
\eqtn
Now the color $c_0(x)$ of  an $m$-tuple $x$ is set to be the tuple $(c'(x^\sigma))_{\sigma\in\mon(M)}$, where $\mon(M)$ is the monoid 
of all mappings $\sigma\colon M\to M$ and  $x^\sigma=(x_{1^\sigma},\ldots,x_{m^\sigma})$.  A linear ordering of 
the colors is inherited in a natural way from 
the lexicographic ordering of the tuples $(1^\sigma,\ldots,m^\sigma)$, $\sigma\in\mon(M)$.

At the second stage, the initial coloring is refined step by step. Namely, if $c_i$ is the coloring constructed at the $i$th step ($i\ge 0$), then the color of an $m$-tuple $x$ in the coloring $c_{i+1}$ is defined to be 
$$
c_{i+1}(x)=(c_i(x),\mltst{(c_i(x_{1\leftarrow \alpha}),\ldots,c_i(x_{m\leftarrow \alpha}))\colon \alpha\in\Omega}),
$$
where $\mltst{\cdot}$ denotes a multiset and, for all $i\in M$, we denote by 
$x_{i\leftarrow \alpha}$ the $m$-tuple $(x_1,\ldots,x_{i-1},\alpha,x_{i+1},\ldots,x_m)$. The algorithm stops when $|\im(c_i)|=|\im(c_{i+1})|$ and the final coloring $c(m,X)$ is set to be~$c_i$.

Two graphs $X$ and $X'$ are said to be {\it $\WL_m$-equivalent} if $\im(c_{X^{}})=\im(c_{X'})$. The \emph{Weisfeiler-Leman dimension} $\dim_{\scriptscriptstyle\WL}(X)$ of a graph $X$ is defined to be the smallest natural $m$ such that every graph $\WL_m$-equivalent to~$X$ is isomorphic to~$X$. 

\subsection{The partition WL$_m$}

The coloring $c_X=c(m,X)$ defines a partition $\WL_m(X)$ of the Cartesian power $\Omega^m$ into the color classes $c_X^{-1}(i)$, where $i$ runs over $\im(c_X)$. This partition is an $m$-ary coherent configuration in the sense of
~\cite{Babai2019}, see also ~\cite{AndresHelfgott2017}. We will not go into details of the general theory of $m$-ary coherent configurations and recall here only a few facts that will be used in the present paper; the interested reader is referred to~\cite[Section~3]{Ponomarenko2022a} or to~\cite[Subsection~3.1]{Chen2023}.

The classes of the partition $\WL_m(X)$ can be considered as $m$-ary relations on~$\Omega$, satisfying certain regularity conditions. These conditions imply the existence of a certain analog of the valency in coherent configurations. More exactly, for a positive integer $k\le m$ and $x\in \Omega^m$, define $\pr_k x:=(x_1,\ldots,x_k)$.
Then, for any class  $\Lambda\in \WL_m(X)$, the number 
\qtnl{250623e}
n_k(\Lambda)=|\{y\in \Lambda\colon \pr_k y=\pr_k x\}|
\eqtn
does not depend on $x\in \Lambda$ (see~\cite[Lemma~3.5]{Chen2023}). 
We extend the definition of $\pr_2$ to subsets $\Lambda\subseteq\Omega^m$ by setting $\pr_2 \Lambda=\{\pr_2 x\colon x\in\Lambda\}$. Then  the $2$-dimensional \emph{projection} 
$\pr_2 \WL_m(X)$ consisting of all sets 
$\pr_2 \Lambda$, $\Lambda\in \WL_m(X)$, is always the set of 
basis relations of a coherent configuration on~$\Omega$. 
Moreover,
\qtnl{250623e1}
\pr_2\WL_2(X)=\WL(X) \qaq \pr_2\WL_k(X)\le \pr_2\WL_m(X)\quad (2\le k\le m).
\eqtn

We complete the subsection by a result \cite[Theorem~3.7]{Ponomarenko2022a} which states a necessary condition 
for the $\WL_m$-equivalence of graphs in terms of 
their coherent configurations.

\thrml{080622a}
Let $X=(\Omega,E)$ and $X'=(\Omega',E')$  be two $\WL_m$-equivalent graphs,  $m\ge 3$, and let  $\cX=\pr_2\WL_m(X)$ and $\cX'=\pr_2\WL_m(X')$. Then for arbitrary points $\alpha_1,\ldots,\alpha_{m-2}\in\Omega$, one can find  some points $\alpha'_1,\ldots,\alpha'_{m-2}\in\Omega'$  and  an algebraic isomorphism 
\qtnl{250623w}
\varphi\colon\cX^{}_{\alpha^{}_1,\ldots,\alpha^{}_{m-2}}\to \cX'_{\alpha'_1,\ldots,\alpha'_{m-2}}
\eqtn
such that $\varphi(1_{\alpha_1^{}})=1_{\alpha_1'},\ \ldots, \varphi(1_{\alpha_{m-2}^{}})=1_{\alpha_{m-2}'}$, and $\varphi(E)=E'$. 
\ethrm

\begin{corollary}\label{coro1}
In the notation of Theorem~\ref{080622a}, 
let $b$ be the base number of $\cX$ and assume that 
$m\geq b+2$ and $m\geq 3$. 
Then $\dim_{\scriptscriptstyle\WL}(X)\le m$.
\end{corollary}
\prf
By assumption, 
one can find some points $\alpha_1,\ldots,\alpha_{m-2}\in\Omega$ 
such that the extension $\cX_{\alpha_1,\ldots,\alpha_{m-2}}$ of 
the coherent configuration $\cX$ with respect to them is discrete. 
By Theorem~\ref{080622a}, one can 
find some points $\alpha'_1,\ldots,\alpha'_{m-2}\in\Omega'$ and algebraic isomorphism~\eqref{250623w} such that $\varphi(E)=E'$. 
Since $\cX_{\alpha_1,\ldots,\alpha_{m-2}}$ is discrete, 
the algebraic isomorphism~$\varphi$ is induced by some  bijection~$f\colon\Omega\to\Omega'$. Therefore, $E^f=E^\varphi=E'$, i.e., $f$ is an isomorphism from $X$ to $X'$. Consequently, any graph $\WL_m$-equivalent to~$X$ is isomorphic to~$X$ and so   $\dim_{\scriptscriptstyle\WL}(X)\le m$.
\eprf

\subsection{Two lemmas}
The two lemmas below are used in the proof of Theorem~\ref{050623a} in Section~\ref{230623j}. The first of them somewhat strengthens the well-known fact that the Weisfeiler-Leman dimension of an arbitrary graph $X$ does not exceed the base number of the coherent configuration $\WL(X)$, increased by~$2$.

\lmml{090623a}
For an integer $k\ge 2$ and a graph $X$, set $b$ to be 
the base number of the coherent configuration $\pr_2 \WL_{k}(X)$. 
Then
$$
\dim_{\scriptscriptstyle\WL}(X)\le\max\{k,b+2\}.
$$
\elmm
\prf
Put $m=\max\{k,b+2\}$. Then the base number of 
$\pr_2\WL_{m}(X)$ is less than or equal to $b$ and 
$m\geq b+2$. 
If $m\geq 3$, then the result follows by 
applying Corollary \ref{coro1} to $\cX=\pr_2\WL_{m}(X)$.
If $m=2$, then $b=0$ and the coherent configuration 
$\WL(X)=\pr_2 \WL_k(X)$ is discrete; thus, it is separable 
and $\dim_{\scriptscriptstyle\WL}(X)\le 2$, see \cite[Theorem~2.5]{Fuhlbr2018a}. 
%
\eprf 

The following lemma could  easily be generalized in different directions (e.g., by replacing two nonadjacent vertices with 
a subgraph of a given type), but we choose the formulation which 
is most relevant to the present paper.

\lmml{090623b}
Let $X$ be a graph  and $e$ an integer. Denote by $s_e(X)$ the relation  of all pairs~$(\alpha, \beta)$ of nonadjacent vertices of $X$, for which $e_X(\alpha,\beta)\ge e$. Then
$$
\pr_2\WL_4(X)\ge \WL(X, s_e(X)).
$$
\elmm
\prf
Denote by $\Delta$ the set of all quadruples $x=(x_1,x_2,x_3,x_4)$ 
of vertices of~$X$, such that $x_1\ne x_2$ and
$$
(x_1,x_2)\not\in E,\quad x_3,x_4\in N_X(x_1,x_2),\quad (x_3,x_4)\in E.
$$ 
It is easily seen that  $\Delta$ is a union of some color classes of the coloring $c_0(m,X)$, see Eq.~\eqref{250623k} for $m=4$. Since every such a class is a union of  some classes of the partition $\fX=\WL_4(X)$, we conclude that $\Delta$ is a union of some classes of~$\fX$.  For each such a class~$\Lambda$, the number 
$$
n_2(\Lambda):=|\{(\alpha,\beta,x_3,x_4)\in\Lambda\colon x_3,x_4\in\Omega\}|
$$
defined by Eq.~\eqref{250623e}, is equal to $e_X(\alpha,\beta)$, and does not depend on the pair $(\alpha,\beta)\in\pr_2 \Lambda$.  Denote by $\Delta^{(e)}$  the union of all classes $\Lambda$ for which $n_2(\Lambda)\ge e$. Then, obviously,
$$
s_e(X)=\pr_2\Delta^{(e)}\in (\pr_2\fX)^\cup.
$$
Since also  $\pr_2\fX\ge \WL(X)$, see \eqref{250623e1}, we conclude that
$$
\pr_2\WL_4(X)=\pr_2\fX\ge \WL(\WL(X),s_e(X))= \WL(X, s_e(X)),
$$
as required.
\eprf

\section{Fon-Der-Flaass graphs from affine planes }

\subsection{The construction}\label{130623h}
The construction described below is a special case of the original Construction~$1$ in \cite{FonDF2002}. In contrast to that paper we  (a) use the affine planes rather than general affine designs, and (b) all these planes are equal to the same plane. 

Let $A$ be an affine plane of order $q$. Denote by $V$ the point set of~$A$  and put $I=\{1,\ldots,q+2\}$. Assume that  for any two distinct indices  $i,j\in I$, we are given a parallel class $\cL_{ij}=\cL_{ji}$ of lines in $A$ and a bijection  $\sigma_{ij}\colon\cL_{ij}\to \cL_{ji}$ such that
\qtnl{140623a}
\cL_{ij}\ne \cL_{ik}\text{ for } k\ne j \qaq \sigma_{ij}^{}=\sigma_{ji}^{-1}.
\eqtn
The first condition enables us to define a $(q+2)\times (q+2)$ array $\cL=(\cL_{ij})$ that can be treated  as a  symmetric Latin square with constant diagonal, in which the off-diagonal elements  are the parallel classes of~$A$. In what follows, we set $\Sigma=\{\sigma_{ij}\}$.

Let us define a graph $X=X_A(\cL,\Sigma)$ with vertex set $\Omega=V\times I$ in which vertices $(v,i)$ and $(u,j)$ are adjacent if and only if $i\ne j$  and $\sigma_{ij}(\bar v)=\bar u$, where $\bar u\in\cL_{ij}$ and $\bar v\in\cL_{ji}$ are the lines of the $A$, containing the points $u$ and $v$, respectively. It was proved in~ \cite{FonDF2002} that the graph $X$ is  strongly regular with parameters \eqref{220623a}. The set of all graphs $X_A(\cL,\Sigma)$ with fixed $A$ and $\cL$, obtained by varying over all possible 
bijections satisfying Eq. \eqref{140623a}, is denoted by $\fF_A(\cL)$; 
when the plane~$A$ is Desarguesian, we use the notation $\fF_q(\cL)$.

The vertex set of any graph $X\in \fF_A(\cL)$ is obviously the disjoint union of the subsets $\Delta_i=V\times\{i\}$ called the \emph{fibers} of $X$; the set of all the $\Delta_i$ is denoted by~$F(X)$. Clearly, $|F(X)|=q+2$. In what follows the points of any fiber are naturally treated as the points  of the plane~$A$.

\lmml{060623a1}
Let $X\in \fF_A(\cL)$ and $F=F(X)$. Then for any two fibers $\Delta,\Gamma\in F$ and any two distinct vertices $\alpha\in\Delta$, $\beta\in\Gamma$, the following statements hold:
\nmrt
\tm{1} if $\Delta\ne\Gamma$ and $\Lambda\in F$, then $|N_X(\alpha,\beta)\cap \Lambda|=0$  or $1$ depending on whether or not $\Lambda\in\{\Delta,\Gamma\}$,
\tm{2} if $\Delta=\Gamma$, then $\alpha$ and $\beta$ are not adjacent and $N_X(\alpha,\beta)\subseteq \Lambda$ for a unique  $\Lambda\in F\setminus\{\Delta,\Gamma\}$.
\enmrt
\elmm
\prf
In what follows we assume that $\Delta=V\times\{i\}$ and  $\Gamma=V\times\{j\}$, where $i,j\in I$, and $\alpha=(v,i)$ and $\beta=(u,j)$ for some $u,v\in V$.

(1) Assume that $i\ne j$. By the construction of $X$, any two distinct points of the same fiber are not adjacent. Hence if $\Lambda\in\{\Delta,\Gamma\}$, then $|N_X(\alpha,\beta)\cap \Lambda|=0$. Now let $\Lambda\notin\{\Delta,\Gamma\}$, i.e., $\Lambda=V\times\{k\}$ for some $k\ne i,j$.  Then 
$$
N_X(\alpha,\beta)\cap \Lambda\subseteq 
\left(\sigma_{ik}(\ov u)\cap \sigma_{jk}(\ov v)\right)\, \times\,\{k\},
$$
where  $\ov u\in \cL_{ik}$ and $\ov v\in \cL_{jk}$  are the lines containing~$u$ and~$v$, respectively. Since $i\ne j$, the lines  $\sigma_{ik}(\ov u)$ and $\sigma_{jk}(\ov v)$ belong to different parallel classes $\cL_{ki}$ and $\cL_{kj}$, and hence have 
exactly one common point. Thus,
$|N_X(\alpha,\beta)\cap \Lambda|=1$.

(2) Let $i=j$. Then the vertices  $\alpha$ and $\beta$ are not adjacent. Let $\gamma=(w,k)$ be a common neighbor of~$\alpha$ and $\beta$. Then the line $\ov u\in\cL_{ik}$ containing~$u$ coincides with the line $\ov v\in \cL_{ik}$ containing~$v$. It follows that $N_X(\alpha,\beta)$ contains all $q$ points $(w',k)$, where $w'$ runs over the points of 
the line $\sigma_{ik}(\ov u)=\ov w\in\cL_{ki}$ containing~$w$. 
Since $n_X(\alpha,\beta)=q$ by Eq. \eqref{220623a}, 
this shows that 
$$
N_X(\alpha,\beta)=\{(w',k)\colon w'\in\ov w\}\subseteq V\times\{k\},
$$
and we are done with $\Lambda=V\times\{k\}$.
\eprf

We complete the subsection by a statement that is used in the proof of the main theorem. It seems that this statement is true not only for the 
graphs in $\fF_A(\cL)$, but also for some other Fon-Der-Flaass 
graphs constructed from affine planes.

\thrml{050623c}
Let $X\in \fF_A(\cL)$ and $\cX\ge \WL(X)$ be a coherent configuration. Assume that  $F(X)\subseteq F(\cX)^\cup$. Then $\cX_{\alpha,\beta}=\cD_\Omega$ for any two  nonadjacent vertices $\alpha$ and $\beta$ of~$X$, not belonging to the same fiber of~$X$. In particular, $b(\cX)\le 2$.
\ethrm
\prf
Let $F(X)=\{\Delta_i\colon i\in I\}$, where $I$ as above. By the assumption of the theorem, each $\Delta_i$ is a homogeneity set of 
the coherent configuration~$\cX$. Since the edge set $E$ of the graph 
$X$ is a relation of the coherent configuration $\WL(X)\le\cX$, 
it follows that if  two indices $i,k\in I$ are distinct, 
then 
$$
s_{ki}=E_{\Delta_k,\Delta_i}\cdot E_{\Delta_i,\Delta_k}
$$ 
is also a relation of  $\cX$. On the other hand, the definition of the graph $X$ implies that~$s_{ki}$ is an equivalence relation on~$\Delta_k$, the classes of which are the lines of the parallel class $\cL_{ki}$ of the affine plane~$A$. \medskip

{\bf Claim. }{\it  Let $T_k$ be the set of all relations $s_{ki}\setminus 1_{\Delta_k}$, $i\ne k$, together with the relation~$1_{\Delta_k}$.  Then  $\cY_k=(\Delta_k,T_k)$ is  the affine scheme associated with 
the plane~$A$.  Moreover, $(\cY_k)_{\alpha,\alpha'}=\cD_{\Delta_k}$ 
for any two distinct points $\alpha,\alpha'\in \Delta_k$.}
\prf
The fact that $\cY_k$ is the affine scheme immediately follows from the first condition in Eq.~\eqref{140623a}. The rest of the statement is a consequence of \cite[Theorem~3.3.8]{CP2019} stating that the extension of an affine scheme with respect to two distinct points is a discrete coherent configuration.
\eprf

Take arbitrary nonadjacent vertices $\alpha\in\Delta_1$ and $\beta\in\Delta_2$. One can see that $\alpha s_{13}=\ov u\times\{1\}$ with $\alpha\in \ov u\in\cL_{13}$, and  $\beta E_{\Delta_2,\Delta_1}=\ov v\times\{1\}$ with some $\ov v\in\cL_{12}$.  By the first condition in Eq.~\eqref{140623a}, the parallel classes $\cL_{13}$ and $\cL_{12}$ are distinct. So  the lines $\ov u$ and $\ov v$ intersect in exactly one point $w$. The vertex $\alpha'=(w,1)$ is different from~$\alpha$, because $\alpha$ and $\beta$ are not adjacent. Moreover,  $\{\alpha'\} = \alpha s_{13} \cap \beta E_{\Delta_2,\Delta_1}$. Since both $\alpha s_{13}$ and $\beta E_{\Delta_2,\Delta_1}$ are homogeneity sets  of the coherent configuration $\cX_{\alpha,\beta}$, this implies that~$\{\alpha'\}$ is a fiber of this configuration, see \cite[Lemma 3.3.5]{CP2019}. In view of the claim, we conclude that
$$
(\cX_{\alpha,\beta})_{\Delta_1}\ge (\cX_{\Delta_1})_{\alpha,\alpha'}\ge(\cY_1)_{\alpha,\alpha'}=\cD_{\Delta_1}.
$$
In a similar way, one can verify that $(\cX_{\alpha,\beta})_{\Delta_2}\ge \cD_{\Delta_2}$.  Consequently, every point of the set $\Delta_1\cup\Delta_2$ forms a fiber of the coherent configuration$~\cX_{\alpha,\beta}$.

To complete the proof it suffices to verify that every point $\lambda\in\Delta$ not belonging to $\Delta_1\cup\Delta_2$ 
 forms a fiber of the coherent configuration $\cX_{\alpha,\beta}$. Indeed, the fiber of it containing~$\lambda$ is a subset of some  $\Delta_i$ for $i>2$.    By the construction of the graph $X$ there are  points $\delta_1\in\Delta_1$ and $\delta_2\in\Delta_2$ such that 
$$
\lambda\in \delta_1 E_{\Delta_1,\Delta_i}\,\cap\,  \delta_2 E_{\Delta_2,\Delta_i}.
$$
Since $\{\delta_1\}$ and $\{\delta_2\}$  are fibers of~ $\cX_{\alpha,\beta}$, the set  on the right-hand side  is a homogeneity one. On the other hand, $\delta_1 E_{\Delta_1,\Delta_i}=\ov w_1\times \{i\}$ and~$\delta_2 E_{\Delta_2,\Delta_i}=\ov w_2\times \{i\}$, where $\ov w_1$ and $\ov w_2$ are distinct nonparallel lines of the affine plane~$A$, see the first condition in Eq.~\eqref{140623a}. Therefore the set $\ov w_1\cap\ov w_2$ and hence the set $\delta_1 E_{\Delta_1,\Delta_i}\cap \delta_2 E_{\Delta_2,\Delta_i}$ is a singleton. Thus, $\{\lambda\}$ is a fiber of~$\cX_{\alpha,\beta}$, as required.
\eprf

\subsection{$4$-condition}
In \cite {Higman1971}, D.~Higman studied the strongly regular graphs of partial geometries.  To formulate one of his results relevant to the present paper, denote by $a(\alpha,\beta)$ the number of the $4$-cliques in a graph~$X$, containing $\alpha$ and $\beta$, and denote by $b(\alpha,\beta)$  the number of diamonds\footnote{The diamond is a complete graph of order~$4$ with one edge removed.} in $X$, in which $\alpha$ and $\beta$ are non-adjacent. 

\lmml{110723a}{\rm \cite[Propositions~6.3 and~6.6(2)]{Higman1971}}
A strongly regular graph with parameters~\eqref{220623a} satisfies the $4$-condition if  and only if $a(\alpha,\beta)=\binom{q}{2}$ for all adjacent vertices $\alpha,\beta$ and $b(\alpha,\beta)=0$ for all nonadjacent vertices $\alpha,\beta$.
\elmm

\crllrl{060623a}
A graph   $X\in\fF_q(\cL)$ satisfies the $4$-condition if the subgraph $X_{N(\alpha,\beta)}$ is complete for any adjacent vertices $\alpha$ and $\beta$ of $X$, where $N(\alpha,\beta)=N_X(\alpha,\beta)$.
\ecrllr
\prf
Let $\alpha$ and $\beta$ be adjacent. Then $a(\alpha,\beta)$ is equal to the edge number $e_X(\alpha,\beta)$ of the subgraph  $X_{N(\alpha,\beta)}$. By the hypothesis, this implies $a(\alpha,\beta)=\binom{q}{2}$. Next, let $\alpha$ and $\beta$ be nonadjacent.  If two distinct vertices $\alpha',\beta'\in N(\alpha,\beta)$ are adjacent, then the graph $X_{N(\alpha',\beta')}$ is not complete (since it contains two nonadjacent vertices $\alpha$ and $\beta$), contrary to the assumption. Therefore the graph  $X_{N(\alpha,\beta)}$ is empty and $b(\alpha,\beta)=0$. Thus the required statement follows from Lemma~\ref{110723a}.
\eprf

\crllrl{060623a14}
If a graph   $X\in\fF_q(\cL)$ satisfies the $4$-condition, then 
$$
e_X(\alpha,\beta)=
\css
0 &\text{if $(\alpha,\beta)\not \in E$,}\\
\binom{q}{2} &\text{if $(\alpha,\beta)\in E$.}\\
\ecss
$$
\ecrllr

\section{Elementary and path switchings}\label{290623a}

\subsection{Elementary switchings}
Throughout this section $A$ is  an affine plane of order~$q$.
  
Let $X\in \fF_A(\cL)$, and let $\Delta=V\times\{i\}$, $\Gamma=V\times\{j\}$ be fibers of~$X$, $i\ne j$.  Let $\cL_{ij}=\cL_{ji}=\{L_1,\ldots,L_q\}$ be a (uniquely determined) parallel class of the plane~$A$.  Given a permutation $f\in\sym(q)$, put 
$$
E'_{\Delta,\Gamma}:=\bigcup_{k=1}^qL_k\times L_{k^f}.
$$
and $E'=(E\setminus E^{}_{\Delta,\Gamma})\cup E'_{\Delta,\Gamma}$. We say that the graph $X'=(\Omega,E')$ is obtained from the graph~$X$ by 
{\it elementary switching} (with respect to $\Delta,\Gamma$) 
if the relation $E^{}_{\Delta,\Gamma}\cup E'_{\Delta,\Gamma}$ is connected, or, equivalently, if the permutation $f$ is a full cycle. 

\lmml{060623w}
Let  a graph $X'$ be obtained from a graph~$X\in\fF_A(\cL)$ by an elementary switching with respect to  $\Delta$ and $\Gamma$. Then $X'\in\fF_A(\cL)$ and, given distinct vertices $\alpha,\beta\in\Omega$, one of the following statements holds:
\nmrt
\tm{1} if $\alpha\ne \beta$ and $\alpha,\beta\not\in \Delta\cup \Gamma$, then $(\alpha,\beta)\in E\ \Leftrightarrow\ (\alpha,\beta)\in E'$, and moreover
$$
N_{X^{}}(\alpha,\beta)=N_{X'}(\alpha,\beta) \qaq |e_{X^{}}(\alpha,\beta)-e_{X'}(\alpha,\beta)|\le 1;
$$
\tm{2}  if $\alpha\in\Delta$  and $\beta\not\in \Delta\cup \Gamma$, then $(\alpha,\beta)\in E\ \Leftrightarrow\ (\alpha,\beta)\in E'$, and moreover
$$
|N_{X^{}}(\alpha,\beta)\setminus N_{X'}(\alpha,\beta)|\le 1 \qaq |e_{X^{}}(\alpha,\beta)-e_{X'}(\alpha,\beta)|\le q-1;
$$
\tm{3}  if $\alpha\in\Delta$  and $\beta\in \Gamma$, then $(\alpha,\beta)\in E\ \Rightarrow\ (\alpha,\beta)\not\in E'$, and moreover
$$
N_{X^{}}(\alpha,\beta)=N_{X'}(\alpha,\beta)  \qaq e_{X^{}}(\alpha,\beta)=e_{X'}(\alpha,\beta);
$$
\tm{4}  if $\alpha,\beta\in\Delta$ or $\alpha,\beta\in\Gamma$,  then $(\alpha,\beta)\not\in E\cup E'$, and moreover,
$$
N_{X^{}}(\alpha,\beta)=N_{X'}(\alpha,\beta) \qoq
N_{X^{}}(\alpha,\beta)\cap  N_{X'}(\alpha,\beta)=\varnothing,
$$
and $e_{X^{}}(\alpha,\beta)=e_{X'}(\alpha,\beta)=0$ 
in either case.
\enmrt
\elmm
\prf 
Obviously, $X'\in \fF_A(\cL)$.

$(1)$ Let $\alpha\ne \beta$ and $\alpha,\beta\not\in \Delta\cup \Gamma$. Since the elementary switching  affects the edges between $\Delta$ and $\Gamma$ only, we have  $(\alpha,\beta)\in E\ \Leftrightarrow\ (\alpha,\beta)\in E'$ and $N_{X^{}}(\alpha,\beta)=N_{X'}(\alpha,\beta)$. Furthermore, if  $\alpha,\beta$ are in the same fiber of $X$ and hence that of $X'$, then  $e_{X^{}}(\alpha,\beta)=e_{X'}(\alpha,\beta)=0$ by Lemma \ref{060623a1}(2). If $\alpha,\beta$ are in different fibers, then 
 $$
 |N_{X^{}}(\alpha,\beta)\cap \Delta|=1 \qaq |N_{X^{}}(\alpha,\beta)\cap \Gamma|=1,
 $$ 
 i.e., there is at most one edge between $N_{X^{}}(\alpha,\beta)\cap \Delta$ and $N_{X^{}}(\alpha,\beta)\cap \Gamma$. Hence, $|e_{X^{}}(\alpha,\beta)-e_{X'}(\alpha,\beta)|\le 1$.

$(2)$ Let $\alpha\in\Delta$  and $\beta\not\in \Delta\cup \Gamma$. Again since the elementary switching affects the edges between $\Delta$ and $\Gamma$ only, we have $(\alpha,\beta)\in E\ \Leftrightarrow\ (\alpha,\beta)\in E'$. By  Lemma~\ref{060623a1}(1),
$$
N_{X}(\alpha,\beta)\cap \Gamma=\{\gamma\}\qaq N_{X'}(\alpha,\beta)\cap \Gamma=\{\gamma '\}
$$
for some vertices $\gamma$ and $\gamma'$, and also 
$$
N_{X^{}}(\alpha,\beta)\setminus \{\gamma\}= N_{X'}(\alpha,\beta)\setminus \{\gamma ' \}
$$ 
which shows that $|N_{X^{}}(\alpha,\beta)\setminus N_{X'}(\alpha,\beta)|\le 1 $. In particular, the graphs induced by $N_{X^{}}(\alpha,\beta)$ in $X$ and  by $N_{X'}(\alpha,\beta)$ in $X'$ differ only in the edges incident to $\gamma$ in~$X$ and~$\gamma'$ in~$X'$. It follows that  $|e_{X^{}}(\alpha,\beta)-e_{X'}(\alpha,\beta)|\le q-1$.

$(3)$ If $\alpha\in\Delta$  and $\beta\in \Gamma$, then the required statements easily follow from  Lemma~\ref{060623a1}(1) and the definition of elementary switching.

$(4)$ Let $\alpha,\beta\in\Delta$. Then  $(\alpha,\beta)\not\in E\cup E'$, and each of the sets $N_{X^{}}(\alpha,\beta)$ and $N_{X^{'}}(\alpha,\beta)$ is contained in a certain fiber, see Lemma \ref{060623a1}(2); in particular, we have $e_{X^{}}(\alpha,\beta)=e_{X'}(\alpha,\beta)=0$.  Without loss of generality, we may assume that the corresponding fibers coincide 
(otherwise, $N_{X}(\alpha,\beta)\,\cap\,N_{X'}(\alpha,\beta)=\varnothing$) and are equal to some $\Lambda\in F(X)$,
$$
N_{X}(\alpha,\beta)\,\cup\,N_{X'}(\alpha,\beta)\subseteq\Lambda.
$$
 Then obviously, $N_{X}(\alpha,\beta)=N_{X'}(\alpha,\beta)$ whenever $\Lambda\ne\Gamma$. Now assume that $\Lambda=\Gamma=V\times\{j\}$.  Then  $\Delta=V\times \{i\}$ for some $i\ne j$, and $\alpha=(u,i)$ and $\beta=(v,i)$ for distinct points $u,v\in V$. It follows that 
 $$
 N_{X^{}}(\alpha,\beta)=L_{k^{}}\times \{j\}\qaq  N_{X'}(\alpha,\beta)=L_{k^f}\times \{j\},
 $$ 
 where $L_{k^{}}\in\cL^{}_{ij}$  is the line of $A$, containing $u$ and $v$, and  $f$ is the permutation in the definition of the elementary switching. Since $L_{k^{}}\cap L_{k^f}=\varnothing$, we conclude that $N_{X^{}}(\alpha,\beta)\cap N_{X^{'}}(\alpha,\beta)=\varnothing$.  
\eprf

\subsection{Path switchings}

A graph $X$ is obtained from a graph $X^*\in\fF_A(\cL)$ 
with edge set $E^*$ by a \emph{path switching}, if there exists a sequence of graphs $X^*=X_0$, $X_1$, $\ldots$, $X_{q+1}=X$ such that $X_i$ is obtained from $X_{i-1}$ by an elementary switching with respect to certain fibers $\Delta_i,\Delta_{i+1}$, $i=1,\ldots,q+1$, and the fibers $\Delta_1,\ldots,\Delta_{q+2}$ are pairwise distinct and 
the same for all graphs $X_i$'s. 

\thrml{080623a}
Let $A$ be an affine plane of even order $q>16$ and $X^*\in\fF_A(\cL)$ a graph satisfying the $4$-condition. Assume that $X$  is a graph obtained from  $X^*$ by a path switching, and set 
\qtnl{050723a}
\cX=\WL(X,s_e(X)),
\eqtn
 where $e=5q-4$ and the relation $s_e(X)$ is defined 
 in Lemma \ref{090623b}. 
 Then  $\cX\ge \WL(X)$ and $F(X)\subseteq F(\cX_\alpha)^\cup$ 
 for some vertex $\alpha$.
\ethrm
\prf
The inclusion $\cX\ge \WL(X)$  is obvious. To prove the second statement, let $\Delta_1,\ldots,\Delta_{q+2}$ and  $X^*=X_0,\ldots, X_{q+1}=X$ be the fibers and the graphs from the definition 
of the path switching. We need two auxiliary lemmas.

\lmml{160623a}
$$
s:=s_e(X)=\bigcup_{k=1}^{q+1}E^*_{\Delta_k,\Delta_{k+1}}.
$$
\elmm
\prf
Let $(\alpha,\beta)\in \Delta_i\times\Delta_j$ for some $i,j\in I$.
Suppose first that the indices $i,j$ satisfy $|i-j|\ne 1$. 
Then for each $k=1,\ldots,q+1$, the pair $(\alpha,\beta)$ satisfies the condition  of one of statements (1), (2), or~(4) of Lemma~\ref{060623w} for $\Delta=\Delta_k$ and $\Gamma=\Delta_{k+1}$. It follows that
$$
e_{X_k}(\alpha,\beta)\le e_{X_{k-1}}(\alpha,\beta)+
\css
1     &\text{if $i\ne j,\  |\{i,j\}\cap\{k,k+1\}|=0$,}\\   
q-1 &\text{if $i\ne j,\  |\{i,j\}\cap\{k,k+1\}|=1$,}\\ 
0    &\text{if $i=j$.}\\
\ecss
$$
The first and second possibilities occur for at most $q-1$ and $4$ values of~$k$, respectively. Consequently, $e_{X^{}}(\alpha,\beta)\le e_{X_0}(\alpha,\beta)+(q-1)+4(q-1)$. Furthermore, if $(\alpha,\beta)\not\in E$, then the assumption  $|i-j|\ne 1$ implies $(\alpha,\beta)\not\in E_0=E^*$. By Corollary~\ref{060623a14}, this yields $e_{X_0}(\alpha,\beta)=0$ whence
$$
e_{X^{}}(\alpha,\beta)\le 5q-5<e,
$$
which shows that $(\alpha,\beta)\not\in s$. In the same way, it 
can also be verified that $(\alpha,\beta)\not\in s$ 
for $|i-j|=1$ under an additional assumption that 
$(\alpha,\beta)\not\in E^*$. 
Thus,
\qtnl{190723a}
s\subseteq \bigcup_{k=1}^{q+1} E^*_{\Delta_k,\Delta_{k+1}}.
\eqtn

To complete the proof, assume that the pair $(\alpha,\beta)$ belongs to the right hand side of \eqref{190723a},  i.e., $j=i+1$.  In this case, this pair can satisfy the condition of any of statements (1)--(3) of Lemma~\ref{060623w} for $\Delta=\Delta_k$ and $\Gamma=\Delta_{k+1}$, where $1\le k\le q+1$. 
By that lemma, we have
$$
e_{X_k}(\alpha,\beta)\ge  e_{X_{k-1}}(\alpha,\beta)-
\css
1     &\text{if $k+1<i$ or $i+1< k$,}\\
q-1 &\text{if $k+1=i$ or $i+1=k$,}\\
0    &\text{if $i=k$.}\\
\ecss
$$
The first and second possibilities occur at most for, respectively, 
$q-1$ and $2$ values of~$k$.
Consequently, $e_{X^{}}(\alpha,\beta)\ge e_{X_0}(\alpha,\beta)-(q-1)-2(q-1)$. On the other hand, since $(\alpha,\beta)\in E^*$, we have  $e_{X_0}(\alpha,\beta)= e_{X^*}(\alpha,\beta)=\binom{q}{2}$ by Corollary~\ref{060623a14}. Since $q>16$, we have
$$
e_{X^{}}(\alpha,\beta)\ge e_{X_0}(\alpha,\beta)-(3q-3)=\binom{q}{2}-(3q-3)>5q-4=e.
$$
This proves that $(\alpha,\beta)\in s$, and hence the inclusion in \eqref{190723a} is an equality, as required.
\eprf

\lmml{160623b}
The set $\Omega_i=\Delta_i\cup\Delta_{q+3-i}$ is a homogeneity set of the coherent configuration $\cX$, $i=1,\ldots,\frac{q+2}{2}$.  Moreover, the subgraph of $(\Omega,s)$ induced by the the set $\Omega_1\cup\Omega_2\cup\Omega_3$ has exactly two connected components and their vertex sets are $\Delta_1\cup \Delta_2\cup\Delta_3$ and $\Delta_q\cup\Delta_{q+1}\cup \Delta_{q+2}$. 
\elmm
\prf
From Lemma \ref{160623a}, it follows that for any point $\alpha\in\Omega$, we have  
\qtnl{190723e}
\alpha s=\css
q  &\text{if $\alpha\in\Omega_1$},\\
2q &\text{otherwise.}\\
\ecss
\eqtn
By Eq. \eqref{180723d} for $d=q$, this implies that $\Omega_1$ is a homogeneity set of the coherent configuration~$\cX$. Assume by  induction that $\Omega_i$ is a homogeneity set of~$\cX$, $i=1,\ldots,k$. Then so is the complement $\Omega'$ of the union $\Omega_1\cup\ldots\cup\Omega_i$ in $\Omega$.   Therefore Eq.~\eqref{190723e} holds for $s$ and  $\Omega_1$ replaced by $s_{\Omega'}$ and $\Omega_{k+1}$, respectively. 
By Eq. \eqref{180723d} for $d=q$, this implies that $\Omega_{k+1}$ is a homogeneity set of~$\cX$, and completes the proof of the first part of the required statement. 

Let us prove the second part of the statement. Since $q>6$, there are no edges between the vertices belonging to $\Delta_1\cup \Delta_2\cup\Delta_3$ and  the vertices belonging to  $\Delta_q\cup\Delta_{q+1}\cup \Delta_{q+2}$. Hence it suffices to verify that the subgraph of the graph $X'=(\Omega,s)$ induced by the the set $\Delta_1\cup \Delta_2\cup\Delta_3$ is connected. Since every vertex of $\Delta_2$ is adjacent in  $X'$ to some vertex of~$\Delta_1$, it suffices to check that any vertices $\alpha\in\Delta_1$ and $\beta\in\Delta_3$ belong to the same component of the graph $X'$.

Let $\alpha=(u,1)$ and $\beta=(v,3)$, where $u$ and $v$ are points of the affine plane~$A$. In view of condition~\eqref{140623a}  the line $\ov u\in\cL_{12}=\cL_{21}$ containing $u$ is different from line $\ov v\in\cL_{32}=\cL_{23}$ containing $v$. So there is a point $w\in \ov u\cap\ov v$. It follows that the vertex $\gamma=(w,2)\in\Delta_2$ is adjacent in $X'$ with both $\alpha$ and $\beta$, as required.
\eprf

From the first statement of Lemma~\ref{160623b}, it follows that $\Lambda=\Omega_1\cup\Omega_2\cup\Omega_3$ is a homogeneity set of the coherent configuration $\cX$. Furthermore, the relation $s_\Lambda$ and hence its  transitive closure $r$  is a relation of~$\cX$. By the second statement of Lemma~\ref{160623b}, we have 
$$
r=(\Delta_1\cup \Delta_2\cup\Delta_3)^2\,\cup\,(\Delta_q\cup\Delta_{q+1}\cup \Delta_{q+2})^2.
$$

Take an arbitrary $\alpha\in \Delta_1$.  Then $\Delta_1=\alpha r\,\cap \Omega_1$ is a homogeneity set of  the coherent configuration $\cX_\alpha$. Assume by induction  that so are the fibers $\Delta_1,\ldots,\Delta_i$ for some $i\le q+1$. By Lemma~\ref{160623a}, we can find the neighborhood of $\Delta_i$ in the graphs $(\Omega,s)$ 
as follows:
$$
\Delta_i s =\Delta_{i-1}\cup\Delta_{i+1}
$$ 
Since $\Delta_i s$ is a homogeneity set of  $\cX_\alpha$, the induction hypothesis implies that so is the set $\Delta_{i+1}=\Delta_i s\setminus \Delta_{i-1}$. 
\eprf

\section{Proof of Theorems~\ref{050623a} }\label{230623j}
The construction of the graph $X^*$ below and Lemma~\ref{080623j} was proposed by M.~Muzychuk (private communication), and, in fact, is a special case of \cite[Proposition~5.1]{Muzychuk2009F} applied  to a hyperoval in a Desarguesian affine plane  of even order.
 
Let $A$ be a Desarguesian affine plane  of even order $q$. Following 
the notation of Subsection~\ref{130623h}, assume that  we are given 
a hyperoval in $A$, say
$$
H=\{h_i\in V\colon i\in I\}.
$$ 
Choose a coordinatization of $A$ so that $H$ does not contains the zero point. For any two distinct indices $i,j\in I$, we put  
\qtnl{220623w}
 \cL_{ij}=\{[h_i-h_j]+v\colon v\in V\},
\eqtn
 where $[h_i-h_j]$  is  the line of $A$, which is incident to zero point and the point $h_i-h_j$.
 
\lmml{150623a}
$\cL_{ij}=\cL_{ji}$ and $\cL_{ij}\ne \cL_{ik}$ for all pairwise distinct indices $i,j,k\in I$.
 \elmm
 \prf
The first statement is obvious.  To prove the second statement, we assume that $\cL_{ij}=\cL_{ik}$ for some pairwise distinct indices $i,j,k\in I$. Then $\cL_{ij}$ and $\cL_{ik}$ have the same line $[u]$ containing zero point; in particular, $[h_i-h_j]=[h_i-h_k]$. Now, the point $h_i$ belongs to the lines  $[u]+h_j$ and $[u]+h_k$ lying in  the same parallel class.  It follows that  they are equal to the same line, and this line contains $h_i$, $h_j$, and~$h_k$, which contradicts the definition of hyperoval.
 \eprf
 
By Lemma~\ref{150623a}, 
the parallel classes~$\cL_{ij}$ and 
the identical bijections $\sigma_{ij}\colon\cL_{ij}\to\cL_{ji}$ 
satisfy Eq.~\eqref{140623a}. 
Hence the array $\cL:=(\cL_{ij})$ 
and the set of bijections $\Sigma:=(\sigma_{ij})$ 
define a strongly regular graph
$$
X^*=X(A,H)=X_A(\cL,\Sigma)
$$
with vertex set $\Omega=V\times I$ and the parameters as in Eq.~\eqref{220623a}, see Subsection~\ref{130623h}.

\lmml{080623j}
The graph $X^*$ belongs to the class $\fF_q(\cL)$, and satisfies the $4$-condition.
\elmm
\prf
The first statement is obvious. To prove the second one, take  arbitrary $a\in \mF_q$ and  $v\in V$, and put
$$
K(a,v)=\{(ah_i+v,i)\in\Omega\colon i\in I\}.
$$

We claim that $K(a,v)$ is a clique (of size $q+2$) of $X^*$. 
Indeed, let $(ah_i+v,i)$ and $(ah_j+v,j)$ be distinct vertices of $K(a,v)$; in particular, $i\ne j$. The difference $(ah_i+v)-(ah_j+v)=a(h_i-h_j)$ belongs to the line $[h_i-h_j]$. 
Consequently, the points $ah_i+v$ and $ah_j+v$ belong to the same  line of the parallel class $\cL_{ij}$, see Eq.~\eqref{220623w}.  But this exactly means that the vertices  $(ah_i+v,i)$ and $(ah_j+v,j)$ are adjacent in $X^*$.

Next, any  two adjacent vertices $\alpha=(v_i,i)$ and~$\beta=(v_j,j)$ of the graph $X^*$ are containd in a common clique $K(a,v)$, where the element $a$ and point $v$ form a solution of the linear system 
$$
\left\{
\begin{aligned}
v_i=ah_i+v,\\
v_j=ah_j+v.\\
\end{aligned}
\right.
$$ 
Hence, $K(a,v)\subseteq N(\alpha,\beta)\cup\{\alpha,\beta\}$, where $N(\alpha,\beta)=N_{X^*}(\alpha,\beta)$. Since the graph~$X^*$ is strongly regular with parameters~\eqref{220623a}, we have $|N(\alpha,\beta) |=q=|K(a,v)|-2$. Thus,
$$
N(\alpha,\beta)\cup\{\alpha,\beta\}=K(a,v),
$$
and hence  the subgraph $N(\alpha,\beta)$ is complete for all adjacent vertices $\alpha$ and $\beta$. By Corollary~\ref{060623a}, this proves that the graph $X^*$ satisfies the $4$-condition.
\eprf

Denote by $\fF=\fF(A,H)$ the family of all graphs $X$ obtained from the graph $X^*=X(A,H)$ by a path switching; in particular, each $X$ is a strongly regular graph with parameters~\eqref {220623a}. 

To prove the first inequality in Eq.~\eqref{220623g}, we need to estimate the Weisfeiler-Leman dimension of an arbitrary graph $X\in\fF$. To this end, let $\cX$  
be the coherent configuration from Theorem~\ref{080623a}, see Eq.~\eqref{050723a}.  The graph $X^*$  satisfies the $4$-condition (Lemma~\ref{080623j}). Consequently, $F(X)\subseteq F(\cX_\alpha)^\cup$  for some vertex~$\alpha$ (Theorem~\ref{080623a}). By Theorem~\ref{050623c} applied to the coherent configuration $\cX_\alpha $, we obtain 
\qtnl{090623c}
\cX_{\alpha,\beta}=(\cX_\alpha)_{\alpha,\beta}=\cD_\Omega
\eqtn
for a suitable vertex~$\beta$. Furthermore, $\pr_2 \WL_4(X)\ge \cX$ 
by Lemma~\ref{090623b}. Together with equality~\eqref{090623c}, 
this shows that
$$
\cD_\Omega\ge (\pr_2 \WL_4(X))_{\alpha,\beta}\ge\cX_{\alpha,\beta}=\cD_\Omega.
$$
It follows that the base number $b$ of the coherent configuration $\pr_2 \WL_4(X)$ is at most~$2$. By Lemma~\ref{090623a}, we finally get 
$$ 
\dim_{\scriptscriptstyle\WL}(X)\le\max\{4,b+2\}=4.
$$

Let us estimate the number of graphs in the family $\fF$. There are $(q+2)!$ ways to choose the sequence of fibers  $\Delta_1,\ldots,\Delta_{q+2}$ for the path switching. As soon as this sequence is fixed, there are $(q-1)!^{q+1}$  ways to choose  $q+1$ full cycles from $\sym(q)$ for the elementary switchings with respect to 
the fibers $\Delta_i$ and $\Delta_{i+1}$. 
Thus the number of distinct graphs in $\fF$ is
$$
(q+2)!\cdot (q-1)!^{q+1}=q^{\mathsf{\Omega}(q^2)}.
$$

To complete the proof, it suffices to verify that every graph $X\in\fF$ 
is isomorphic to at most $n^2=\mathsf{\Omega}(q^6)$ graphs 
from~$\fF$. Indeed, 
Eq.~\eqref{090623c} shows that every isomorphism $\pi$ from $X$ to 
an arbitrary graph on~$\Omega$ is uniquely determined by the choice of 
the vertices $\alpha^\pi$ and $\beta^\pi$: in fact, if $\pi'$ is 
another isomorphism such that $\alpha^\pi=\alpha^{\pi'}$ and $\beta^\pi=\beta^{\pi'}$, then
$$
\pi'\,\pi^{-1}\in \aut(\WL(X,1_\alpha,1_\beta))=\aut(\cX_{\alpha,\beta})=\aut(\cD_\Omega)=\{\id_\Omega\},
$$
i.e., $\pi=\pi'$. It follows that among all the graphs on~$\Omega$, and hence among all graphs in~$\fF$, there are at most $|\Omega|^2$ graphs isomorphic to~$X$.

\subsection*{Acknowledgements}
The research of Jinzhuan Cai, Jin Guo, and Ilia Ponomarenko 
is supported by the National Natural Science Foundation of 
China (Grant No. 11961017).
The research of Alexander Gavrilyuk is supported by JSPS KAKENHI Grant Number 22K03403.

\providecommand{\bysame}{\leavevmode\hbox to3em{\hrulefill}\thinspace}
\providecommand{\MR}{\relax\ifhmode\unskip\space\fi MR }
\providecommand{\MRhref}[2]{%
	\href{http://www.ams.org/mathscinet-getitem?mr=#1}{#2}
}
\providecommand{\href}[2]{#2}

\bibliographystyle{amsplain}
\bibliography{wffsrg}

\end{document}